\documentclass[12pt]{article}
\usepackage[english]{babel}
\usepackage{amssymb,amsmath,graphics}
\parskip\medskipamount
\newtheorem{theorem}{Theorem}[section]
\newtheorem{Theorem}{Main Result}
\newtheorem{remark}[theorem]{Remark}

\newtheorem{lemma}[theorem]{Lemma}
\newtheorem{cor}[theorem]{Corollary}
\def\<{\langle}
\def\>{\rangle}
\newcommand{\proof}{\emph{Proof.~}}

\newcommand{\cF}{\mathcal{F}}
\newcommand{\dd}{\mathsf{d}}

\def\qed{{\hfill\hphantom{.}\nobreak\hfill$\Box$}}

\newcommand{\A}{\mathbf{A}}
\newcommand{\R}{\mathbb{R}}

\newcommand{\N}{\mathbb{N}}
\newcommand{\F}{\mathbb{F}}
\newcommand{\inc}{\mbox{\tt I}}

\begin{document}

\author{Koen Struyve\thanks{The first author is supported by  the Fund for Scientific Research --
Flanders (FWO - Vlaanderen)} \and Hendrik Van Maldeghem\thanks{The second author is partly supported by a
Research Grant of the Fund for Scientific Research -- Flanders (FWO - Vlaanderen)} }
\title{\bf Generalized polygons with non-discrete valuation defined by two-dimensional affine $\R$-buildings}

%\date{}
\maketitle

\begin{abstract}
\noindent
In this paper, we show that the building at infinity of a two-dimensional affine $\R$-building is a generalized
polygon endowed with a valuation satisfying some specific axioms. Specializing to the discrete case of affine
buildings, this solves part of a long standing conjecture about affine buildings of type
$\widetilde{\mathsf{G}}_2$, and it reproves the results obtained mainly by the second author for types
$\widetilde{\mathsf{A}}_2$ and $\widetilde{\mathsf{C}}_2$. The techniques are completely different from the
ones employed in the discrete case, but they are considerably shorter, and general (i.e., independent of the
type of the two-dimensional $\R$-building).
\end{abstract}

\section{Introduction}
Buildings, introduced by Jacques Tits in the early 60's, are geometrical structures that help us understand
large classes of groups, such as semisimple algebraic groups, semisimple Lie groups, Chevalley groups and their
twisted analogues, groups of mixed type, classical groups, Kac-Moody groups, arithmetic groups, etc.
Fundamental theorems in this interaction are the classification of spherical buildings of rank at least 3
\cite{Tit:74} and of affine buildings of rank at least 4 \cite{tits84}, both by Jacques Tits, the
classification of 2-spherical twin buildings of rank at least 3 for which the rank 2 residues are not
isomorphic to the ``smallest spherical buildings of rank 2'' by Bernhard M\"uhlherr \cite{Mue:02}, and the
classification of Moufang spherical buildings of rank 2 by Jacques Tits and Richard Weiss \cite{Tit-Wei:02}.
From geometric point of view, the buildings that only just fail the hypotheses of these theorems are important
as well for several reasons. First of all, their existence proves that the assumptions of the classification
theorems are necessary. Secondly, the structure of these exceptions is usually less restricted and hence it
tells us in greater generality how buildings look like. Moreover, these buildings have geometric invariants and
characteristic properties of their own which make them useful for a variety of other geometric or group
theoretic applications. A standard example of this behaviour are the non-classical projective planes (which
are spherical buildings of rank 2), in particular the translation planes and the associated spreads in
projective space.

The subject of the present paper is related to affine buildings of rank 3. However, the original statement of
the above mentioned classification theorem of affine buildings by Jacques Tits is slightly more general in that
it also includes structures that very much look like affine buildings, but are not. These structures are called
\emph{affine apartment systems} (``syst\`emes d' appartements'') by Jacques Tits, but later on people also
called them \emph{non-discrete affine buildings} \cite{ronan} or \emph{$\mathbb{R}$-buildings}. They recently
received a renewed attention because they provide examples of CAT(0)-spaces.

In a series of rather long papers \cite{Mal:87,Mal:88,Mal:89,Mal:90,Han-Mal:90}, the second author (jointly
with Guy Hanssens in the last quoted paper) investigates in detail two classes of rank 3 affine buildings
(namely, those of type $\widetilde{\mathsf{A}}_2$ and those of type $\widetilde{\mathsf{C}}_2$) and
characterizes the corresponding spherical buildings at infinity (which are projective planes and generalized
quadrangles, respectively). This leads to many new examples of such affine buildings, explicitly defined and
with knowledge of the automorphism groups. Originally, the characterization made use of the notion of a
\emph{discrete valuation} on the algebraic structures that coordinatize projective planes and generalized
quadrangles, but in later papers \cite{Mal:92,Mal:93}, the valuation was defined directly on the geometry. The
hope was that with such a direct definition, the case of type $\widetilde{ \mathsf{G}}_2$, which was the only
remaining case, would become treatable with much less effort. One of the reasons why it did \emph{not} is that,
although the paper \cite{Mal:92} provides the exact condition for a generalized hexagon with valuation, the
lack of symmetry in the formulae prevents us from deducing a \emph{general} formulae independent of the type,
and hence from (1) further generalization to non-discrete valuations, and (2) composing a type-free proof.

In this paper, we start such a type-free approach, which ought to eventually lead to a characterization of  all
irreducible 2-dimensional affine apartment systems.   In fact, we present half of such a characterization
here. More in particular, we will show how any irreducible 2-dimensional affine apartment system gives
rise to a generalized polygon with a specific valuation, by which we mean, with the terminology of
\cite{Mal:92}, an explicitly defined weight sequence. One of the crucial observations to achieve this is to
slightly modify, or re-scale, the valuation as defined from a rank 3 affine building as defined in
\cite{Mal:93}. Indeed, roughly speaking, the valuation between two elements as defined in \cite{Mal:93} counted the graph
theoretic distance between two vertices in the simplicial complex related to the affine building. The purpose
was to end up with a natural number. But taking the Euclidean distance instead will put much more symmetry into
the picture, and at the same time we will have a closed formula for the weight sequences. Also the non-discrete
case can clearly be included in a natural way. The fact that the discrete case enjoys a characterization as in
\cite{Mal:92} seems to be a happy coincidence in this viewpoint.

\section{Preliminaries and Main Results}\label{MR}
\subsection{$\R$-buildings}
Let $(\overline{W},S)$ be a finite irreducible Coxeter system. So $\overline{W}$ is presented by the set $S$ of
involutions subject to the relations which specify the order of the products of every pair of involutions. This
group has a natural action on a vector space $V$ of dimension $|S|$. Let $\A$ be the affine space associated to
$V$. We define $W$ to be the group generated by the translations of $\mathbf{A}$ and $\overline{W}$.

Let $\mathcal{H}_0$ be the set of hyperplanes of $V$ corresponding to the axes of the reflections in $S$ and
all its conjugates. Let $\mathcal{H}$ be the set of all translates of all elements of $\mathcal{H}_0$. The
elements of $\mathcal{H}$ are called \emph{walls} and the (closed) half spaces they bound are called
\emph{half-apartments} or \emph{roots}. A \emph{vector sector} is the intersection of all roots that (1) are
bounded by elements of $\mathcal{H}_0$, and (2) contain a given point $x$ that does not belong to any element
of $\mathcal{H}_0$. The bounding walls of these roots will be referred to as the \emph{side-walls} of the
vector sector. A vector sector can also be defined as the closure of a connected component of
$V\setminus(\cup\mathcal{H}_0)$. Any translate of a vector sector is a \emph{sector}, with corresponding
translated \emph{side-walls}. A \emph{sector-facet} is an infinite intersection of a given sector with a finite
number of its side-walls. This number can be zero, in which case the sector-facet is the sector itself; if this
number is one, then we call the sector-facet a \emph{sector-panel}. The intersection of a sector with all its
side-walls is a point which is called the \emph{source} of the sector, and of every sector-facet defined from
it.

An \emph{$\R$-building} (also called an \emph{affine apartment system}) (definition by Jacques Tits as can be
found in~\cite{ronan} by Mark Ronan, along with some historic background) is an object $(\Lambda,\cF)$
consisting of a set $\Lambda$ together with a collection $\cF$ of injections of $\mathbf{A}$ into $\Lambda$
obeying the five conditions below. The image of $\mathbf{A}$ under a $f\in \cF$ will be called an
\emph{apartment}, and the image of a sector, half-apartment, \dots of $\mathbf{A}$ under a certain $f\in \cF$
will be called a \emph{sector}, \emph{half-apartment}, \dots of $\Lambda$.
\begin{itemize}
 \item[(A1)] If $w\in W$ and $f\in \cF$, then $f \circ w \in \cF$.
 \item[(A2)] If $f,f' \in \cF$, then $X=f^{-1}(f'(\mathbf{A}))$ is closed and convex
 %the intersection of a finite number of roots
 in $\textbf{A}$, and $f|_X = f'\circ w|_X$ for some $w\in W$.
 \item[(A3)] Any two points of $\Lambda$ lie in a common apartment.
\end{itemize}
The last two axioms allow us to define a function $\dd : \Lambda \times \Lambda \rightarrow \R^+$ such that for
any $a,b \in \A$ and $f\in \cF$, $\dd(f(a),f(b))$ is equal to the Euclidean distance between $a$ and $b$ in
$\A$.
\begin{itemize}
 \item[(A4)] Any two sectors contain subsectors lying in a common apartment.
 \item[(A5')] Given $f \in \cF$ and a point $\alpha \in \Lambda$, there is a retraction
 $\rho : \Lambda \rightarrow f(\A)$ such that the preimage of $\alpha$ is $\{\alpha\}$ and which diminishes $\dd$.
\end{itemize}

We call $|S|$, which is also equal to $\dim \mathbf{A}$, the \emph{dimension} of $(\Lambda,\cF)$. We will
usually denote $(\Lambda,\cF)$ briefly by $\Lambda$, with slight abuse of notation.

A detailed analysis of this definition and variations of it has carried out by Anne Parreau in~\cite{parreau}.
In particular, she shows that $\dd$ is a distance function.

One can associate spherical buildings of type $(\overline{W},S)$ to these $\R$-buildings in two ways. The first
way to do so is to construct the building at infinity. Two sector-facets of $\Lambda$ will be called
\emph{asymptotic} if the distance between them is bounded. Due to the triangular inequality this is an equivalence
relation. The equivalence classes (named \emph{facets at infinity}) form a spherical building $\Lambda_\infty$
of type $(\overline{W},S)$ called the \emph{building at infinity} of $(\Lambda,\cF)$. The chambers of
$\Lambda_\infty$ are the equivalence classes of asymptotic sectors. An apartment $\Sigma$ of $\Lambda$
corresponds to an apartment $\Sigma_\infty$ of $\Lambda_\infty$ in a bijective relation. 

A second way to construct a spherical building is to look at the `local' structure instead of the asymptotic
one. Let $\alpha$ be a point of $\Lambda$, and $F,F'$ two sector-facets with source $\alpha$. Then these two
facets will \emph{locally coincide} if their intersection is a neighbourhood of $\alpha$ in both $F$ and $F'$.
This relation forms an equivalence relation defining \emph{germs of facets} as equivalence classes (notation
$[F]_\alpha$). These germs form a (possibly non-thick) building  $[\Lambda]_\alpha$ of type $(\overline{W},S)$,
called the \emph{residue} at $\alpha$. If $\Sigma$ is an apartment containing $\alpha$, then $[\Sigma]_\alpha$
will be used to denote the corresponding apartment in $[\Lambda]_\alpha$.

The following lemma by Anne Parreau will prove to be an important tool in our proofs.
 \begin{lemma}[Parreau~\cite{parreau}, Proposition 1.8]\label{lemma:parr}
Let $x$ be a chamber of the building at infinity $\Lambda_\infty$ and $C$ a sector with source $\alpha \in
\Lambda$. Then there exists an apartment $\Sigma$ containing an element of the germ $[C]_\alpha$ and such that
$\Sigma_\infty$ contains $x$. \qed
 \end{lemma}

This has also an interesting corollary.
 \begin{cor}[Parreau~\cite{parreau}, Corollary 1.9]
Let $\alpha$ be a point of $\Lambda$ and $F_\infty$ a facet of the building at infinity. Then there is a unique
facet $F' \in F_\infty$ with source $\alpha$. \qed
 \end{cor}

The unique facet of the previous corollary will be denoted as $(F_\infty)_\alpha$ or $F_\alpha$.

If the dimension of $\Lambda$ is at least 3, then $\Lambda_\infty$ is a spherical Moufang building and, in
principle, $\Lambda$ is known, see \cite{tits84}. Hence, in this paper we will only deal with $\R$-buildings of
dimension 2, i.e., $|S|=2$ and $\overline{W}$ is the dihedral group of order $2n$, for some $n\in\mathbb{N}$,
$n\ge 3$. So the building at infinity and the residues are generalized $n$-gons, in which we can discern two
types of sector-panels, one called \emph{points} $P$, the others \emph{lines} $L$ and an \emph{incidence
relation} $\inc$ between both (the choice which type of sector-panels are the points or lines can be chosen
arbitrarily). An \emph{element} of a generalized polygon will just be any point or line. \emph{Collinearity}
and \emph{concurrency} are defined in the usual way, \emph{adjacency} is both of these relations barring
equality. Roman letters will be used for elements of the building at infinity, Greek letters for points of
$\Lambda$.

Let $x,y$ be two adjacent elements of $\Lambda_\infty$ and $\alpha \in \Lambda$, then we denote the length (measured with the
distance $\dd$) of the common part of the sector-panels $x_\alpha$ and $y_\alpha$ by $u_\alpha(x,y)$. The
mapping $u_\alpha$, defined on adjacent pairs of elements of $\Lambda_\infty$, has some properties very similar
to that of a valuation in the sense of \cite{Mal:93}. After we will have introduced these in
subsection~\ref{valuation}, we will be able to state our main result.

\subsection{Generalized polygons}

Generalized polygons are the geometries corresponding to the spherical rank 2 buildings. Since we will use some
specific terminology of these geometries, we introduce this now.

A \emph{generalized $n$-gon}, $n\in\mathbb{N}$, $n\ge 2$, or \emph{generalized polygon} $\Gamma=(P,L,\inc)$ is
a structure consisting of a \emph{point set} $P$, a \emph{line set} $L$ (with $P\cap L=\emptyset$), and a
symmetric \emph{incidence relation} $\inc$ between $P$ and $L$, turning $P\cup L$ into a bipartite graph
$\mathfrak{G}$ satisfying the following axioms.
\begin{itemize}
 \item[(GP1)] Every element is incident with at least two other elements.
 \item[(GP2)] For every pair of elements $x,y\in P\cup L$,
 there exists a sequence $x_0=x,x_1,\ldots,x_{k-1},x_k=y$, with $x_{i-1}\inc x_{i}$ for $1\leq i\leq k$ and
 with $k\leq n$.
 \item[(GP3)] The sequence in (GP2) is unique whenever $k<n$.
\end{itemize}

Axiom (GP3) implies that the smallest positive integer $k$ for which there exists a sequence $x_0\inc
x_1\inc\cdots\inc x_k\inc x_0$ of $k$ different points and lines is $2n$. In this case, the subgeometry induced
by these elements is an \emph{apartment} of the generalized $n$-gon. If not all elements are distinct, but
still $k=2n$, then we will sometimes call the geometry induced by these elements a \emph{degenerate apartment}.

\subsection{Generalized polygons with (non-discrete) valuation}\label{valuation}

Let $\Gamma=(P,L,\inc)$ be a generalized $n$-gon with point set $P$ and line set $L$,  and let $u$ be a
function called a \emph{valuation} acting on both pairs of collinear points and pairs of concurrent lines, and
images in $\R^+ \cup \{ \infty \} $ (we use the natural order on this set with $\infty$ as largest element).
Then we call $(\Gamma,u)$ an \emph{$n$-gon with (non-discrete) valuation} with \emph{weight sequence}
$(a_1,a_2, \dots, a_{n-1}, a_{n+1},a_{n+2}, \dots a_{2n-1} ) \in (\R^+)^{2n-2}$ if the following conditions are
met :
\begin{itemize}
 \item[(U1)] On each line there exists a pair of points $p$ and $q$ such that $u(p,q)=0$ and dually for
points.
 \item[(U2)] $u(x,y) = \infty$ if and only if $x=y$.
 \item[(U3)] $u(x,y) < u(y,z)$ implies $u(x,z) = u(x,y)$ if $x,y$ and $z$ are collinear points or concurrent lines.
 \item[(U4)] Whenever $x_1 \inc x_2 \inc\dots \inc  x_{2n} \inc x_1$, with $x_i \in P \cup L$, one has
\begin{equation}
\sum_{i=1}^{n-1} a_i u(x_{i-1},x_{i+1}) = \sum_{i=n+1}^{2n-1} a_i u(x_{i-1},x_{i+1}).
\end{equation}
\end{itemize}
One direct implication of (U3) is that $u$ is symmetric (by putting $x=z$). Also remark that this definition is
self-dual when interchanging lines and points, so whenever a statement is proven, we also have proven the dual
statement. Finally, we note that, due to (U2), axiom (U4) is trivially satisfied whenever the $x_i$, $1\leq
i\leq 2n$, form a degenerate apartment.

The difference with the definition in \cite{Mal:92} is that in the current definition, the element $x_1$ is
arbitrary, while in \cite{Mal:92}, $x_1$ was required to be a point. On the other hand, in \cite{Mal:92}, the
image of $u$ had to be the natural numbers together with $\infty$. The main result of \cite{Mal:92} says that,
in this case, $n\in\{3,4,6\}$, the function $u$ is also a valuation on the dual $n$-gon, and the
weight-sequences are uniquely determined up to duality. These weight-sequences are, however, only self-dual if
$n=3$. Hence, only in the case $n=3$, a valuation on an $n$-gon in the sense of \cite{Mal:92} will be a
valuation on an $n$-gon in the above sense. However, rescaling the valuation between lines by a factor
$\sqrt{2}$ (multiplying or dividing according to the weight-sequence) for $n=4$ turns the valuation on a
$4$-gon in the sense of \cite{Mal:92} into a valuation on in the above sense. Similarly for $6$-gons. Taking
this rescaling into account, we see that the above definition is essentially a generalization of the definition
in \cite{Mal:92}. We will come back to this in more detail in Section~\ref{discrete}, where we will show how
our main result proves one direction of the conjectures stated in \cite{Mal:92} and \cite{Mal:93}.

Our main result reads:

\subsection{Main Result}
\begin{Theorem}
Let $(\Lambda,\cF)$ be an $\R$-building and $\alpha \in \Lambda$. Then $u_\alpha$ defines a valuation on
$\Lambda_\infty$ with weight sequence $(a_1,a_2, \dots, a_{n-1}, a_{n+1},a_{n+2}, \dots a_{2n-1} )$, where $a_i
= |\sin(i \pi / n) |$.
\end{Theorem}

\section{Proof of the main result}
The first lemma deals with the exact shape of the intersection of two sectors with same source and sharing a
sector-facet.

 \begin{lemma}\label{lemma:intersect}
Let $C$ and $C'$ be two sectors with the same source $\tau$ which share a sector-facet $F$. Then the
intersection of both is formed by the convex hull of $F$ and the common part of the other two sector-facets of
$C$ and $C'$.
 \end{lemma}
\proof Take any apartment $\Sigma$ containing $C$ (and thus also $\tau$). If $\Sigma$ contains $C'$ then there
is nothing left to prove. If this is not the case then there is a unique apartment $\Sigma'_\infty$ at infinity
containing $C'$ and sharing a half-apartment with $\Sigma_\infty$. A remark in~\cite{parreau} (page 10) states
that if two apartments share a half-apartment at infinity, they also do in the building itself. This implies
the exact form of the intersection. \qed

If $C\neq C'$, then such an intersection is called a \emph{chimney} with source $\tau$ by Guy Rousseau
(\cite{rousseau}). We refer to the \emph{width} of the chimney as the distance between the parallel walls
bordering it.
 \begin{cor}\label{cor:width}
Let $r,s,t$ be elements of $\Lambda_\infty$ such that $r \inc s \inc t$, and $\tau$ a point, then the width of
the chimney defined as the intersection of the sector containing $r_\tau$ and $s_\tau$, and the one containing
$s_\tau$ and $t_\tau$, is equal to $ \sin(\pi/n) u_\tau(r,t)$.
 \end{cor}
\proof Directly from the definitions and the previous lemma. \qed

Now let $\alpha$ be an arbitrary point of $\Lambda$ and consider the mapping $u_\alpha$.  The axiom (U1) will
be satisfied because given an element $x$ at infinity there is always an apartment containing $x_\alpha$ where
we can find the needed element $y$ such that $u_\alpha(x,y)=0$. The second axiom (U2) follows trivially and
(U3) follows from the convexity of sector-panels.

The main difficulty is (U4).  Let $x_0$ and $x_n$ be two opposite elements of $\Lambda_\infty$ and $M:=(x_0,
x_1, \dots ,x_n) \in (P \cup B)^{n+1}$  such that $x_0 \inc x_1 \inc \dots \inc x_n$. % For each $l \in \R^+$,
We define the function $$f:\R^+\rightarrow\R^+:l\mapsto\sum_{i=1}^{n-1} \sin(i \pi / n)
u_\beta(x_{i-1},x_{i+1}),$$ with $\beta \in (x_0)_\alpha$ at distance $l$ from $\alpha$. If we can prove that
$f$ only depends on $x_0,x_n$ and $\alpha$, then we have proven (U4) and the main result (in view of the fact
that (U4) is trivially in degenerate apartments).

%The next lemma implies that $f$ is continuous.
%\begin{lemma}
%Let $x,y,z$ be elements of $\Lambda_\infty$, and let $x$ and $y$ be adjacent. Then the function $z_\alpha
%\rightarrow \R^+ : \beta \mapsto u_\beta(x,y)$ is continuous.
%\end{lemma}
%\proof Let $\beta$ be on $z_\alpha$. Due to lemma~\ref{lemma:parr} we can find an apartment $\Sigma$ containing
%$x_\beta$ and an element of $[z]_\beta$. Similarly, one has an apartment $\Sigma'$ containing $y_\beta$ and an
%element of $[z]_\beta$. The intersection of both of these apartments contains the common part of $x_\beta$ and
%$y_\beta$, and nothing more of either $x_\beta$ or $y_\beta$, and it also contains an interval $I_\beta$ of
%$z_\beta$ containing $\beta$.

%Let $\gamma$ be a point in $I_\beta$. The sector panels $x_\gamma$ and $y_\gamma$ lie in the apartments
%$\Sigma$ and $\Sigma'$ respectively, because of this $x_\gamma \cap y_\gamma$ will be equal to $x_\gamma \cap
%\Sigma' \cap \Sigma$. As $\Sigma \cap \Sigma'$ is closed and convex in $\Sigma$ the function will be convex,
%and thus continuous (Proposition 1.19 in~\cite{phelps})\qed

% Because of the previous lemma, one can say that $|u_\beta(x,y) - u_\gamma(x,y)| \leq d(\beta,\gamma) / \tan(\pi /n) $. This implies continuity of our function in one direction, the other direction is proven analougously. \qed

Before we go on we need the notion of ``distance in the residues''. Let $x$ and $y$ be elements of
$\Lambda_\infty$ and $\beta \in \Lambda$. Then we define the \emph{residual distance $\dd_\beta(x,y)$ at $\beta$} to be
the distance between $[x]_\beta$ and $[y]_\beta$ as defined in the generalized $n$-gon $[\Lambda]_\beta$ (a
point and an incident line are at distance 1, two collinear points are at distance 2, \dots).

The next lemma investigates the local behaviour of the valuations.

 \begin{lemma}\label{lemma:progress}
Let $r,s,t$ be elements of $\Lambda_\infty$ such that $r \inc s \inc t$, and $\beta$ a point on $(x_0)_\alpha$
with $\dd(\alpha,\beta) =l$. Then there exists $\delta > 0$ such that for any $\beta'$ on $(x_0)_\alpha$ with
$\dd(\alpha,\beta') \in [l,l+\delta]$, the following holds~: $$u_{\beta'}(r,t) = u_\beta(r,t) + \epsilon\frac{
\sin(\dd_\beta(s,x_0 ) \pi /n)  }{ \sin(\pi/n)} \dd(\beta,\beta'),$$ where $\epsilon$ is a constant equal to
 $$\begin{cases}
 -1 \mbox{ if }\dd_\beta(r,x_0 ) = \dd_\beta(t,x_0 ) =\dd_\beta(s,x_0 )-1, \\
 1 \mbox{ if } \dd_\beta(r,x_0 ) = \dd_\beta(t,x_0 )=\dd_\beta(s,x_0 )+1, \\
 0 \mbox{ if } \dd_\beta(r,x_0 ) \neq \dd_\beta(t,x_0 ) .
\end{cases}$$
\end{lemma}
\proof Let $C$ be the sector spanned by $r_\beta$ and $s_\beta$ and $C'$ the one by $s_\beta$ and $t_\beta$.
Both these sectors have source $\beta$. Using lemma~\ref{lemma:parr}, we can find apartments $\Sigma$ and
$\Sigma'$ containing $C$ and an element of the germ $[x_0]_\beta$, and $C'$ and an element of the germ
$[x_0]_\beta$, respectively.  Let $\delta$ be the length of $(x_0)_\beta$ included in $\Sigma\cap\Sigma'$.
Obviously $\delta > 0$. Let $\beta'$ be on $(x_0)_\alpha$ with $\dd(\alpha,\beta') \in [l,l+\delta]$. The
sectors $C_{\beta'}$ and ${C'}_{\beta'} $ with source $\beta'$ now lie in the apartments $\Sigma$ and
$\Sigma'$, respectively.

Using the intersection of both apartments one can easily calculate that the width of the chimney defined by
$r,s$ and $t$ with source at $\beta'$ is $\epsilon \sin(\dd_\beta(s,x_0 ) \pi /n) d(\beta,\beta') )$ larger
than the one with source $\beta$, with $\epsilon$ as in the table above. Using corollary~\ref{cor:width} we now
obtain the desired result. \qed

As an immediate consequence of the previous lemma, we see that $f$ is continuous.

%The included figures now show the apartment $\Sigma$ with the blue areas that are certainly in the intersection, the red certainly not by standard convexity and closedness arguments.  A few %straightforward calculations now lead to $u_{\beta'}(r,t) = u_\beta(r,t) + \epsilon (\sin(d_\beta(s,x_0 ) \pi /n  ) / \sin(\pi/n)) d(\beta,\beta')$ with the correct $\epsilon$. \qed

Applying the previous lemma to the (finite number) of valuations occurring in the definition of $f$ now implies
that for every $l \in \R^+$ there exists $\bar{\delta}>0$ (the minimum occurring in the application to each
valuation) and  $a_l \in \R$ such that $f(l') = f(l) + a_l(l'-l)$ for every $l' \in [l,l+\bar{\delta}]$. The
next step in our proof is to show that $a_l$ only depends on $x_0, x_n, \alpha$ and $l$. One thing that is
directly clear is that $a_l$ only depends on the distances $\dd_\beta(x_0,x_i)$ with $i \in \{1,2, \dots,n \}$,
and on $\beta$ on $(x_0)_\alpha$ via $d(\alpha,\beta)=l$. Because of this we can reduce this combinatorially as
follows. Define the sequence $(y_0,y_1,\ldots,y_n)$, with $y_i:=\dd_\beta(x_0,x_i)$, $i \in \{0,1,2, \dots,n
\}$. This sequence consists of non-negative integers such that two consecutive ones differ by exactly one, and
the extremities $y_0$, which equals $0$, and $y_n$ are constants. An entry different from the extremities with
the property that both neighbours are strictly smaller will be called a \emph{peak}; if both neighbours are
strictly larger, then we call the entry a \emph{valley}. The sequence will determine the $a_l$ uniquely.

If two sequences produce the same $a_l$ we will say that they are equivalent. We now show that each sequence is equivalent to the unique sequence with no valleys, which will be called the standard sequence. Therefore we look
at the sum $\chi$ of all the $y_i$'s. The number $\chi$ is clearly an integer and bounded. Consider any sequence
different from the standard sequence, then it has at least one valley, say at the entry $y_j=m$. We now break
the problem down to some different cases and show that in each case that the given sequence is equivalent with
one obtained from the first one by replacing $y_j$ by $y_j+2$. This equivalent sequence has larger sum,
and because this sum is an integer and is bounded by the sum obtained from the standard sequence, recursion implies that all sequences are equivalent to the standard sequence. Note that $j\geq 2$, so $j-2$ is always
well-defined.

In the following we will denote $\pi/n$ by $\pi_n$ for ease of notation.

\begin{itemize}

\item[$(i)$] \textbf{Case $(y_{j-2},y_{j-1},y_j,y_{j+1},y_{j+2})=(m+2,m+1,m,m+1,m)$.} So, as indicated above,
we show that this is equivalent with $(y_{j-2},y_{j-1},y'_j,y_{j+1},y_{j+2})=(m+2,m+1,m+2,m+1,m)$. Indeed,
using the expression for $a_l$ from the definition of $f$ and Lemma~\ref{lemma:progress}, we see that we must
show \begin{eqnarray*} - \sin({j\pi_n}) \sin({m\pi_n}) +\sin({(j+1)\pi_n})\sin({(m+1)\pi_n})=\\ -
\sin({(j-1)\pi_n})\sin({(m+1)\pi_n}) + \sin({j\pi_n})\sin({(m+2)\pi_n}).\end{eqnarray*}

Indeed, we perform the following elementary calculations.
\begin{align*}
& - \sin(j\pi_n ) \sin(m\pi_n) +sin((j+1)\pi_n ) \sin((m+1)\pi_n)  \\
&= 1/2 ( - \cos((j-m)\pi_n) +\cos((j +m)\pi_n)  +\cos((j-m)\pi_n)  -\cos((j+m+2)\pi_n)  ) \\
&= 1/2 (  \cos((j+m)\pi_n)  -\cos((j+m+2)\pi_n)  ),
\end{align*}
while
\begin{align*}
& - \sin((j-1)\pi_n ) \sin((m+1)\pi_n) +sin(j\pi_n ) \sin((m+2)\pi_n)  \\
&= 1/2 ( - \cos((j-m-2)\pi_n) +\cos((j +m)\pi_n)  +\cos((j-m-2)\pi_n) -\cos((j+m+2)\pi_n)  )   \\
&= 1/2 (  \cos((j+m)\pi_n)  -\cos((j+m+2)\pi_n)  ).
\end{align*}
It follows that the two sequences are equivalent.

\item[$(ii)$] \textbf{Case $(y_{j-2},y_{j-1},y_j,y_{j+1},y_{j+2})=(m,m+1,m,m+1,m+2)$.} This is directly
analogous to the previous case.

\item[$(iii)$] \textbf{Case $(y_{j-2},y_{j-1},y_j,y_{j+1},y_{j+2})=(m,m+1,m,m+1,m)$.} Here, we show that this
is equivalent with $(y_{j-2},y_{j-1},y'_j,y_{j+1},y_{j+2})=(m+2,m+1,m+2,m+1,m)$. Indeed, as before, we must
show that \begin{eqnarray*}\sin((j-1)\pi_n )\sin((m+1)\pi_n) - \sin(j\pi_n ) \sin(m\pi_n) +\sin((j+1)\pi_n )
\sin((m+1)\pi_n)\\ \phantom{\sin((j-1)\pi_n )\sin((m+1)\pi_n)}=\sin(j\pi_n ) \sin((m+2)\pi_n).\end{eqnarray*}
This equality is the same as the one in Case~$(i)$, but with one term swapped from side. The same conclusion
follows.

\item[$(iv)$] \textbf{Case $(y_{j-2},y_{j-1},y_j,y_{j+1},y_{j+2})=(m+2,m+1,m,m+1,m+2)$.} Here we must show that
$$\begin{array}{l}-\sin(j\pi_n ) \sin(m\pi_n)=\\ -\sin((j-1)\pi_n ) \sin((m+1)\pi_n) + \sin(j\pi_n )
\sin((m+2)\pi_n) -\sin((j+1)\pi_n ) \sin((m+1)\pi_n).\end{array}$$ This equality is the same as in Case $(iii)$
but with $m$ substituted by $-m-2$. Again the same conclusion follows.

\item[$(v)$] \textbf{Case $j=n-1$}. In this case we can reuse the previous arguments by adding an extra element
$x_{n+1} \inc x_n$ with corresponding $y_{n+1} := y_{n} \pm 1$, and extending $f$ with an extra coefficient
$\sin(n \pi / n) u_\beta(x_{n-1},x_{n+1})$ (which is zero anyway due to $\sin\pi=0$).
\end{itemize}

This proves that each sequence is equivalent to the standard sequence, and thus that all sequences are equivalent and $a_l$ only depends on $x_0, x_n, \alpha$ and $l$.

We now need an elementary result from analysis, which we prove for completeness' sake.

\begin{lemma} \label{lemma:const}
If $g$ is a continuous real function defined over $\R^+$ such that for every $l \in \R^+$ there is a $\delta$
such that $g(l') = g(l)$ for every $l' \in [l,l+\delta]$, then $g$ is constant over $\R^+$.
\end{lemma}
\proof Define $\Psi := \{ x \in \R^+  |  (\exists \delta'
>0)(\forall x' \in [x-\delta',x+\delta'] )(g(x) =g(x'))  \}$ as the set of ``constant points''. If an interval lies
completely in $\Psi$, then $g$ is constant over that interval because the preimage of the image of an element
in such an interval is both open (due to the definition of constant points) and closed (because of continuity)
in the connected interval. By using continuity, this is also true for the closure of an interval lying
completely in $\Psi$. If the set $\R^+ \backslash \Psi$ is non-empty, then it has an infimum $t$. Note that by
assumption, there exists $\delta>0$ such that $[0,\delta[\subseteq \Psi$. Hence $t>0$ and the interval $[0,t[$
lies completely in $\Psi$, implying that $g$ is constant over $[0,t]$. But we also know that there exists
$\delta'$ such that $g$ is constant over $[t,t+\delta']$, so $[0,t+\delta']$ lies in $\Psi$. This contradicts
the fact that $t$ is an infimum. So $\Psi = \R^+$ and $g$ is constant over $\R^+$. \qed

\begin{lemma}
There is an $l \in \R^+ $ such that $f(l') =0$ if $l' \geq l$.
\end{lemma}
\proof Let $i$ be minimum with respect to the property $\dd_\alpha(x_0,x_i) \neq i$. It is clear that, if
$\beta \in (x_0)_\alpha$, then $\dd_\beta(x_0,x_j) =j $ for $j<i$ (because the sectors spanned by $x_0$ till
$x_j$ with source $\alpha$ form a part of an apartment and contain those with source $\beta$). Suppose there is
no $\beta \in (x_0)_\alpha$ such that also $\dd_\beta(x_0,x_i) =i$. In such a case we have that the function
$$g: \R^+ \rightarrow \R^+ : l \mapsto u_\beta(x_{i-2},x_i), \mbox{ with }\dd(\beta,\alpha)=l,$$ is strictly positive
for each $\beta \in (x_0)_\alpha$ (because a zero value would imply that $\dd_\beta(x_0,x_i) =i$). As we know
by lemma~\ref{lemma:progress}, for every $l\in \R^+$ there is a $\delta$ such that $$g(l') =g(l) -
\frac{\sin((i-1) \pi_n  )}{\sin\pi_n} (l'-l), \mbox{ for every }l' \in [l,l+\delta].$$  The function $g(l) +
\frac{\sin((i-1) \pi _n  )}{\sin\pi_n} l$ then complies to the statement of Lemma~\ref{lemma:const} and is thus
constant. But this is impossible since for $l$ large enough, this would imply that $g(l)$ is negative.
Consequently $g$ cannot be strictly positive, yielding that there is a $\beta \in (x_0)_\alpha$ such that also
$\dd_\beta(x_0,x_i) =i$.

Repeating this process a finite number of times will produce an $l$ such that  $\dd_\beta(x_0,x_n) =n$ if
$\dd(\beta,\alpha) \geq l$. This implies that $u_\beta(x_{i-1},x_{i+1})$ is zero for each $i \in \{1,2, \dots,
n-1\}$, which on its turn implies that $f(d(\beta,\alpha))=0$. \qed

\begin{remark}\em The above lemma can also be proved (in a possibly simpler way) using some additional notions introduced in \cite{tits84},
in particular, the tree $T(W)$ corresponding to the ``wall'' $W=\{x_0,x_n\}$ in $\Lambda_\infty$, but we have
chosen the proof above since it avoids additional definitions. Briefly, one knows that the sector-panel
$(x_0)_\alpha$ intersects some element of the $\R$-tree $T(W)$, say $\beta$ is in this intersection. Then
obviously we can put $l=\dd(\alpha,\beta)$.
\end{remark}

Let us reiterate what we know about the function $f$ defined over $\R^+$:
\begin{itemize}
\item[(O)] For high enough values it is zero.%
\item[(C)] The function is continuous.%
\item[(P)] For every $l \in \R^+$ there is a $\bar{\delta}$ and an $a_l \in \R$ such that $f(l') = f(l) +
a_l(l'-l)$ for every $l' \in [l,l+\bar{\delta}]$ where $a_l$ depends only on $l, x_0,x_n$ and $\alpha$.
\end{itemize}
\begin{lemma}
Two functions satisfying the three conditions {\rm (O), (C)} and {\rm (P)} (with the same $a_l$) are equal over
$\R^+$.
\end{lemma}
\proof Because we know that $f$ satisfies the above conditions, we can assume that one of the functions is $f$,
let the other be $f'$. Consider $g=f'-f$, then $g$ is continuous, is for high enough values zero, and for every
$l \in \R^+$ there is a $\delta$ (the minimum of the two $\bar{\delta}$ related to $f$ and $f'$) such that
$g(l') = g(l)$ for every $l' \in [l,l+\delta]$. Lemma~\ref{lemma:const} now implies that $g$ is constant, and
thus zero over $\R^+$.

%We now prove that such a function $g$ is equal to zero over $\R^+$. Define $\Psi := \{ x \in \R^+ |  (\exists \delta' >0)(\forall x' \in [x-\delta,x+\delta] )(g(x) =g(x'))  \}$ as the set of constant points. If an interval lies completely in $\Psi$, then $g$ is constant over that interval because the preimage in the interval of an element is both open (due to the definition of constant points) and closed (because of continuity) in the connected interval. By using continuity, this is also true for the closure of an interval lying completely in $\Psi$. If the set $\R^+ \backslash \Psi$ is non-empty, then it has an infinum $t$. In that case the interval $[0,t[$ lies completely in $\Psi$, implying that $g$ is constant over $[0,t]$. But we also know that there is an $\delta$ such that $g$ is constant over $[t,t+\delta]$, so $[0,t+\delta]$ lies in $\Psi$ contradicting the fact that $t$ is an infinum. So $\Psi = \R^+$ and $g$ is constant, and thus zero over $\R^+$.

This implies that $f$ and $f'$ are equal. \qed

As $a_l$ only depends on $l, x_0,x_n$ and $\alpha$, the previous lemma has the corollary that $f$ only depends
on  $x_0,x_n$ and $\alpha$, which has previously been said to imply (U4). This completes the proof of our Main
Result.

\section{The discrete case}\label{discrete}

Let (U4$'$) be the Condition (U4) with the additional requirement that $x_0\in P$, and let $\Gamma=(P,L,\inc)$
be a generalized $n$-gon, $n\ge 3$. Suppose that $(\Gamma,u)$ satisfies (U1), (U2), (U3) and (U4$'$), and
suppose in addition that the image of $u$ is $\N \cup \{\infty\}$, the set of natural numbers, including $0$, together with the element $\infty$. Then we say that $(\Gamma,u)$ is a generalized polygon with discrete valuation. The main result of \cite{Mal:92} says that, in
this case,  $n\in\{3,4,6\}$ and the weight sequence $(a_1,a_2, \dots, a_{n-1}, a_{n+1},a_{n+2}, \dots a_{2n-1}
)$ can be chosen as follows.
\begin{itemize}
\item[(WS3)] If $n=3$, then $(a_1,a_2,a_4,a_5)=(1,1,1,1)$.%
\item[(WS4)] If $n=4$, then $(a_1,a_2,a_3,a_5,a_6,a_7)\in\{(1,1,1,1,1,1),(1,2,1,1,2,1)\}$.%
\item[(WS6)] If $n=6$, then
$(a_1,a_2,\ldots,a_5,a_7,\ldots,a_{11})\in\{(1,1,2,1,1,1,1,2,1,1),(1,3,2,3,1,1,3,2,3,1)\}$.
\end{itemize}

In the cases (WS4) and (WS6), where there are two possibilities, it is proved in \cite{Mal:92} that these
weight sequences are dual to one another, i.e., if $(\Gamma,u)$ has one weight sequence, then, if $\Gamma^*$ is
the dual of $\Gamma$ (obtained from $\Gamma$ by interchanging the point set with the line set), then
$(\Gamma^*,u)$ is a polygon with discrete valuation with respect to the other weight sequence.

This gave birth to the conjecture that \emph{a generalized hexagon $\Gamma$ is ``isomorphic'' to the building
at infinity of some (thick) affine building of type $\widetilde{\mathsf{G}}_2$ if and only if there exists $u$
such that $(\Gamma,u)$ is a generalized hexagon with discrete valuation and one of the two above weight
sequences}. The Main Result of the present paper seems to be in contradiction with this, since, applied to
discrete affine buildings of type $\widetilde{\mathsf{G}}_2$, there is only one weight sequence, namely
$$(a_1,a_2,\ldots,a_5,a_7,\ldots,a_{11})=\left(\frac{1}{2},\frac{\sqrt{3}}{2},1,\frac{\sqrt{3}}{2},
\frac{1}{2},\frac{1}{2},\frac{\sqrt{3}}{2},1,\frac{\sqrt{3}}{2},\frac{1}{2}\right),$$ and it does not consist
of only natural numbers! But the above conjecture was evidenced by the situation for the types
$\widetilde{\mathsf{A}}_2$ and $\widetilde{\mathsf{C}}_2$, where the valuation measured simplicial distance,
and not Euclidean distance, as in the present paper. In the $\widetilde{\mathsf{G}}_2$ case, this means that,
in view of the fact that the lengths of the panels (of a chamber) containing the special vertex (for
terminology, see \cite{tits84}) have ratio $1:\sqrt{3}$, to go from the weight sequence of the present paper to
the weight sequences of the discrete valuation, we must multiply the valuation on the point pairs with
$\sqrt{3}$ (or do this with the valuation on line pairs), and then take a suitable multiple.

As explained earlier, one can do a similar procedure with type $\widetilde{\mathsf{C}}_2$, as is clear from the
above.

\section{Algebraization and an example}

For the types $\widetilde{\mathsf{A}}_2$ and $\widetilde{\mathsf{C}}_2$, it is shown in
\cite{Mal:87,Mal:88,Mal:89,Mal:90,Mal:93} that a generalized $n$-gon with discrete valuation, $n=3,4$, is
equivalent with the corresponding general coordinatizing structure with discrete valuation. These
coordinatizing structures are ``planar ternary rings'' and ``quadratic quaternary rings''. The proofs of these
equivalences can be taken over verbatim for the current situation of non-discrete valuations. This enables us
to give some non-classical examples of generalized $n$-gons with non-discrete valuation in the current sense,
and for reasons of simplicity, we restrict ourselves to the case $n=3$.

EXAMPLE. Let $\F$ be a field, and denote by $\F\<t\>$ the integral domain containing all formal polynomials
$\sum_{j\in J}a_jt^j$, with $J$ some finite subset of $\R^+$ (depending on the element of $\F\<t\>$). Denote by
$\F\{t\}$ the corresponding quotient field. Let $\sigma$ be an automorphism of $\F$ of finite order and denote
by the same symbol the automorphism of $\F\{t\}$ induced by $\sigma$ in the natural way, i.e.,
$(\sum_{j\in J}a_jt^j )^\sigma=\sum_{j\in J}a_j^\sigma t^j$. Then $|\<\sigma\>|$ is finite. With
these data, and somewhat similar to Example~7.2 in \cite{Mal:87}, we can construct an Andr\'e quasifield,
which, with the obvious valuation, is a quasifield with non-discrete valuation, and hence gives rise to a
projective plane with non-discrete valuation which is not a Moufang plane. The details of this construction are
left to the reader. These examples will only become really interesting as soon as the converse of the Main
Result is proven to be true, i.e., as soon as we know that projective planes with non-discrete valuations give
rise to affine $\R$-buildings whose buildings at infinity are precisely these projective planes. This will be
carried out in a future paper, where we will also come back in some more detail to the current examples.

\end{document}